\documentclass[10pt,twoside]{classe-ma_e}

\usepackage{amscd}      
\usepackage{amssymb}     
\usepackage[intlimits]{amsmath}     
\usepackage{xypic}      
\LaTeXdiagrams          
\usepackage[all,v2]{xy}      
\xyoption{2cell}
\UseAllTwocells
\xyoption{frame}
\CompileMatrices

            {\nolinebreak $\Box$ \end{trivlist}}

\newcommand{\noprint}[1]{}

\renewcommand{\tilde}{\widetilde}

\newcommand{\toto}{\rightrightarrows}

\newcommand{\upst}{^{\ast}}

\newcommand{\com}{^{\scriptscriptstyle\bullet}}
\newcommand{\lcom}{_{\scriptscriptstyle\bullet}}

\newcommand{\XX}{{\mathfrak X}}

\newcommand{\RR}{{\mathfrak R}}

\newcommand{\Gg}{{\mathfrak g}}

\newcommand{\rr}{{\mathbb R}}

\newcommand{\del}{\partial}

\newcommand{\ldiag}[1]%
       {\makebox[0cm]{${\scriptstyle#1}\downarrow\phantom{\scriptstyle#1}$}}
\newcommand{\ldiagup}[1]%
       {\makebox[0cm]{${\scriptstyle#1}\uparrow\phantom{\scriptstyle#1}$}}
\newcommand{\rdiag}[1]%
       {\makebox[0cm]{$\phantom{\scriptstyle#1}\downarrow{\scriptstyle#1}$}}
\newcommand{\sediagr}[1]%
       {\makebox[0cm]{$\phantom{\scriptstyle#1}\searrow{\scriptstyle#1}$}}
\newcommand{\nediagr}[1]%
       {\makebox[0cm]{$\phantom{\scriptstyle#1}\nearrow{\scriptstyle#1}$}}
\newcommand{\rdiagup}[1]%
       {\makebox[0cm]{$\phantom{\scriptstyle#1}\uparrow{\scriptstyle#1}$}}
\newcommand{\swdiag}[1]%
       {\makebox[0cm]{$\phantom{\scriptstyle#1}\swarrow{\scriptstyle#1}$}}
\newcommand{\sediag}[1]%
       {\makebox[0cm]{${\scriptstyle#1}\searrow\phantom{\scriptstyle#1}$}}
\newcommand{\nediag}[1]%
       {\makebox[0cm]{${\scriptstyle#1}\nearrow\phantom{\scriptstyle#1}$}}

\newcommand{\longiso}{\stackrel{\textstyle\sim}{\longrightarrow}}

\newcommand{\doublearrowstack}[2]%
                      {{{{\scriptstyle#1}\atop{\textstyle\longrightarrow}}\atop{{\textstyle\longrightarrow}\atop{\scriptstyle#2}}}}
\newcommand{\rightleftarrowstack}[2]%
                      {{{{\scriptstyle#1}\atop{\textstyle\longrightarrow}}\atop{{\textstyle\longleftarrow}\atop{\scriptstyle#2}}}}
\newcommand{\leftrightarrowstack}[2]%
                      {{{{\scriptstyle#1}\atop{\textstyle\longleftarrow}}\atop{{\textstyle\longrightarrow}\atop{\scriptstyle#2}}}}

\newcommand{\overtoparrow}%
{\makebox[0cm]{\beginpicture
\setcoordinatesystem units <.8cm,.4cm> point at 0 0
\setplotarea x from -3 to 3, y from 0 to 1
\setquadratic
\plot -3 0 0 1 3 0 /
\put{\vector(3,-1){0}}[Bl] at 3 0
\endpicture}}

\newcommand{\underbottomarrow}%
{\makebox[0cm]{\beginpicture
\setcoordinatesystem units <.8cm,.4cm> point at 0 0
\setplotarea x from -3 to 3, y from 0 to 1
\setquadratic
\plot -3 1 0 0 3 1 /
\put{\vector(3,1){0}}[Bl] at 3 1
\endpicture}}

\newcommand{\ses}[5]%
{0\longrightarrow#1\stackrel{#2}{ \longrightarrow}#3\stackrel{#4}{
\longrightarrow}#5\longrightarrow0}

\newcommand{\dt}[6]%
{#1\stackrel{#2}{longrightarrow}#3 \stackrel{#4}{\longrightarrow}#5
\stackrel{#6}{\longrightarrow} #1[1]}  
 
\newcommand{\cat}[1]%
{(\mbox{\rm #1})}

\newcommand{\reals}{{\Bbb      R}}

\newcommand{\gm}{\Gamma}
\newcommand{\lon }{\longrightarrow }

\newcommand{\btheta}{\bar{\theta}}

\newcommand{\frakg}{\Gg}
\newcommand{\ra}{\rangle}
\newcommand{\la}{\langle}
\newcommand{\half}{\frac{1}{2}}
\newcommand{\hol}{\mbox{Hol}}    

\def\gpd{\,\lower1pt\hbox{$\longrightarrow$}\hskip-.24in\raise2pt
             \hbox{$\longrightarrow$}\,}

\newtheorem{them}{Theorem}[section]

\newtheorem{prop}[them]{Proposition}
\newtheorem{cor}[them]{Corollary}

\newtheorem{numrmk}{\it Remark\/}

\newtheorem{theoreme}{Th\'eor\`eme}[section]

\ComParit{S0764-4442}
\PIT{FLA}
\PXMA{????}
\Add{?}
\Volume{332}
\Year{2001}
\FirstPage{1}
\LastPage{??}
\AuteurCourant{Running authors}
\TitreCourant{Running title} 
\Journal
\Rubrique{Rubrique}{Heading}
\SousRubrique{Sous-rubrique}{Sub-Heading}
\PresentePar{First name}{NAME}
\Recu{jour mois ann\'ee}{apr\`es r\'evision}{jour mois ann\'ee}
\begin{document}
\title{%
Equivariant Gerbes over Compact Simple  Lie Groups
}
\author{%
Kai Behrend~$^{\text{a}}$,\ \
Ping Xu~$^{\text{b}}$,\ \
Bin Zhang~$^{\text{c}}$
}
\address{%
\begin{itemize}\labelsep=2mm\leftskip=-5mm
\item[$^{\text{a}}$]
University of British Columbia \\
E-mail: behrend@math.ubc.ca
\item[$^{\text{b}}$]
Pennsylvania State University \\
E-mail: ping@math.psu.edu 
\item[$^{\text{c}}$]
State University of New York at Stony Brook \\
E-mail: bzhang@math.sunysb.edu
\end{itemize}
}
\maketitle
\thispagestyle{empty}
\begin{Abstract}{%
Using groupoid $S^1$-central extensions, we present,
for a compact simple  Lie group $G$,
 an infinite dimensional model of $S^1$-gerbe
over the differential stack $G/G$ whose Dixmier-Douady
class corresponds to the canonical
 generator of the equivariant cohomology $H_G^3 (G)$.
}\end{Abstract}
\begin{Ftitle}{%
Gerbes Equivariantes sur les groupes de Lie simples compacts
}\end{Ftitle}
\begin{Resume}{%
En utilisant des extensions $S^1$-centrales de groupo\"{\i}des, nous pr\'esentons,
dans le cas d'un groupe simple compact $G$ un mod\`ele de dimension infinie
d'une $S^1$-gerbe sur un champ diff\'erentiable $G/G$ dont la
classe de
Dixmier-Douady
correspond au g\'en\'erateur canonique de la cohomologie \'equivariante $H_G^3 (G)$.                      
}\end{Resume}

\AFv

Soit $G$ un groupe de Lie compact simple, le groupe de cohomologie
\'equivariante $H^3_G (G)$ contient un g\'en\'erateur canonique
$[\omega +\Omega] $ dont la classe est enti\`ere (voir
Sec. \ref{sec:AMM}  pour la d\'efinition
 de $\omega$ et $\Omega$),
o\`u $G$ agit sur lui-m\^eme par conjugaison.
Nous r\'ealisons cette classe en termes d'une extension $S^1$-centrale
de groupo\"{\i}des, ou en tant que la classe de Dixmier-Douady d'une gerbe
sur un champ diff\'erentiable. Le champ est celui correspondant \`a la
transformation de groupo\"{\i}de $G\times G\toto G$.
Cet exemple s'int\`egre dans la th\'eorie g\'en\'erale d\'evelopp\'ee
dans \cite{BX, BX2}. Notre construction est divis\'ee en deux \'etapes.
La premi\`ere \'etape s'inscrit dans le cadre de la g\'eom\'etrie de
Poisson: pour une vari\'et\'e de Poisson affine $\frakg^*$ induite
par un 2-cocycle d'alg\`ebre de Lie $\lambda \in \wedge^2 \frakg^*$,
nous construisons son groupo\"{\i}de symplectique ainsi qu'une
extension $S^1$-centrale
de groupo\"{\i}des. La construction se fait de la mani\`ere suivante.
Soit $S^{1}\lon \tilde{G}\stackrel{\pi}{\lon} G$ une extension
$S^1$-centrale au niveau des groupes de Lie. On forme les groupo\"{\i}des
de transformations $\gm (=G\times \frakg^*  ) \toto \frakg^* $
et $R (=\tilde{G}\times \frakg^* ) \toto \frakg^* $
($G$ agit sur $\frakg^*$ par l'action de jauge \eqref{eq:coad2},
tandis que $\tilde{G}$ agit sur $\frakg^*$ essentiellement par la m\^eme
action obtenue par composition avec le morphisme de groupes
$\pi : \tilde{G}\to G$). En consid\'erant  $R$ comme sous-groupo\"{\i}de
du groupo\"{\i}de symplectique $T^* \tilde{G}\toto \widetilde{\frakg}^*$,
on peut consid\'erer les images r\'eciproques sur $R$ des formes symplectique
et
de Liouville sur $T^* \tilde{G}$ et obtenir repectivement
une 2-forme ferm\'ee $\omega_R \in \Omega^2 (R) $ et une un-forme
$\theta_R \in \Omega^1 (R)$. Soit $\tilde{\pi}: R\to \gm$ la projection
naturelle induite par $\pi : \tilde{G}\to G$.

\begin{theoreme}
\begin{enumerate}
\item La 2-forme ferm\'ee  $\omega_R $ est basique relativement
au $S^1$-fibr\'e
   $R\to \gm$, et donc se projette en une 2-forme ferm\'ee
 $\omega_{\gm}$ sur $\gm$, c.\`a.d,
 $\omega_{R}=\tilde{\pi}^*\omega_{\gm}$.
\item $\omega_{\gm}$ est
 symplectique, ce qui est compatible avec la structure de
 groupo\"{\i}de de telle sorte que cela d\'efinisse une structure de
 groupo\"{\i}de symplectique sur $\gm$.  C'est le groupo\"{\i}de
 symplectique de la vari\'et\'e de Poisson affine $\frakg^*$.
\item
 $\tilde{\pi} : R\lon \gm $ est une extension $S^1$-centrale de
 groupo\"{\i}des de Lie.  \item $\theta_R$ est une forme de connexion
 de pr\'equantification sur le $S^1$-fibr\'e $\tilde{\pi } : R\to \gm $
 compatible avec la structure de groupo\"{\i}de, c.\`a.d., $\partial
 \theta_R =0$ et $d \theta_R = \tilde{\pi}^* \omega_{\gm}$.
\end{enumerate}
\end{theoreme}

En d'autres termes, $\theta_R\in C^2_{DR}(R\lcom ) $ 
est une
pseudo-connexion sur $R$  au sens de
\cite{BX} dont $\omega_{\gm}\in Z^3_{DR}(\gm\lcom )$
 est la  pseudo-courbure.

Pour un groupe de Lie simple compact $G$ d'alg\`ebre de Lie
$\frakg$, la forme basique  sur $\frakg$ induit une 2-cocycle naturel
d'alg\`ebre de Lie sur l'alg\`ebre de Lie de lacets $L\frakg$.
En appliquant la construction ci-dessus, on obtient un groupo\"{\i}de
symplectique $({LG} \times L\frakg \toto L\frakg,
\omega_{LG \times L\frakg} ) $ ainsi qu'une extension $S^1$-centrale
$\widetilde{LG} \times L\frakg \toto L\frakg $ de ce dernier.
Ici nous identifions $L\frakg$ avec  $L\frakg^*$ par la forme de Killing.
L'application d'holonomie $\hol : L\frakg\to G$ induit un morphisme naturel
de Morita du groupo\"{\i}de ${LG} \times L\frakg \toto L\frakg $
vers $G\times G\toto G$,   dont le morphisme induit en cohomologie de de Rham
envoie  $[\omega +\Omega] $
sur la classe  $[\omega_{LG \times L\frakg}]$.
Nous prouvons donc que :

\begin{theoreme}
Soit $\omega +\Omega$ un 3-cocycle comme ci-dessus, d\'efinissant une classe
enti\`ere dans $H^3_{G}(G)$.
L'extension  $S^1$-centrale  $\widetilde{LG} \times L\frakg \toto L\frakg $
de ${LG} \times L\frakg \toto L\frakg $
correspond  \`a une  $S^1$-gerbe sur le champ $G/G $ de classe de
Dixmier-Douady $[\omega +\Omega ]\in H^3_{G}(G)$.
\end{theoreme}

\section{Introduction}
Let $G$ be a compact simple Lie group, the equivariant
cohomology  group $H^3_G (G)$ contains a canonical generator
of integer
 class, where $G$ acts on itself
by conjugation. We realize this class in terms of  a
Lie groupoid $S^1$-central extension, or as the Dixmier-Douady
class of an $S^1$-gerbe over a differential stack.
The stack is the one  which corresponds to the transformation
groupoid $G\times G\toto G$.  This example fits into the general 
theory developed in \cite{BX, BX2}. Recall that associated
to every Lie groupoid  $\Gamma\toto M$, there are De Rham cohomology
groups defined as follows.
 Define $\Gamma_p=\underbrace{\Gamma\times_M\ldots\times_M\Gamma}_{\text{$p$ times}}$, i.e., $\Gamma_p$ is the manifold of composable
 sequences of $p$ arrows in the groupoid $\Gamma\toto M$
 ($\Gamma _1=\Gamma, \Gamma _0=M$).
 We have $p+1$ canonical maps $\gm_p\to \gm_{p-1}$ (each leaving out one of the $p+1$ objects involved a sequence of composable arrows), giving rise to a         diagram
$\xymatrix{
\ldots \gm_2
\ar[r]\ar@<1ex>[r]\ar@<-1ex>[r] & \gm_1\ar@<-.5ex>[r]\ar@<.5ex>[r]
&\gm_0\,.}$
In fact,  $\gm \lcom$ is a simplicial manifold.   We introduce the double complex
$\Omega\com(\Gamma\lcom)$. Its boundary maps are $d:
\Omega^{k}( \gm_p ) \to \Omega^{k+1}( \gm_p )$, the usual exterior
derivative of differential forms and $\partial
:\Omega^{k}( \gm_p ) \to \Omega^{k}( \gm_{p+1} )$,  the alternating
sum of the pull back maps of the above diagram.
We denote the total complex by $C^*_{DR}(\Gamma\lcom)$ and
the total differential by $\delta=(-1)^pd+\del$.
The total cohomology groups of $\Omega\com(\Gamma\lcom)$,
$H_{DR}^k(\Gamma\lcom)=H^k\big(\Omega\com(\Gamma\lcom)\big)$
are called the {\em De~Rham cohomology }groups of $\Gamma\toto M$.          
In the case that $\Gamma\toto M$ is a transformation groupoid
$G\times M\toto M$, these are the $G$-equivariant cohomology groups. 
In \cite{BX}, we discussed the general
 question of how to realize a De Rham integer 3-cocycle
  in terms of  an analogue of curvature
on a  Lie groupoid $S^1$-central extension.
When this  3-cocycle consists of only one term
$\omega \in \Omega^2 (\gm )$, the 3-cocycle 
condition  is equivalent to that $(\gm\toto M, \omega )$ 
is  a symplectic groupoid if $\omega $ is further assumed
to be non-degenerate. $S^1$-central extensions of a symplectic groupoid
were studied extensively by Weinstein and one of us in \cite{WX},
to which we refer the readers for details.

When $\gm \toto M$ is the transformation groupoid  $G\times M\toto M$,
then $\omega +\Omega \in C_{DR}^3 (\gm\lcom )$, where
  $\omega \in \Omega^2 (\gm )$ and
$\Omega\in \Omega^3 (M)$,  is a 3-cocycle if and only if 
 $d\Omega=0,\   d\omega =\alpha^* \Omega -\beta^* \Omega \ \mbox{ and }
 (\partial^*_{0}- \partial^*_{1}+\partial^*_{2})\omega=0$. 
In this case, the class of $\omega +\Omega $
defines an element in the   
$G$-equivariant cohomology  group $H^3_{G}(M)$. When $M=G$  is a compact
simple Lie group and $G$ acts
on itself by conjugation,
 an explicit formula for both an $\omega \in \Omega^2 (\gm )$
and an $\Omega\in \Omega^3 (M)$ appeared in \cite{AMM} in the
study of group valued momentum maps and the moduli spaces
of flat connections over two surfaces.
 See also \cite{Weinstein}
and \cite{JHGW} for the related topics and motivation. 
 However, 
 the fact that $\omega +\Omega$  is
a 3-cocycle in $ C_{DR}^3 (\gm\lcom )$ was overlooked  in the
literature.
In this Note, we reinstate this fact, and show
that $[\omega +\Omega]$ is an integer class by constructing
an  $S^1$-gerbe over the stack $G/G$ (also called
a $G$-equivariant gerbe over $G$)  
which has $[\omega +\Omega]$ as its Dixmier-Douady class.
Our method is to pass from  the groupoid $G\times G\toto G$
to a Morita equivalent 
infinite dimensional symplectic  groupoid, where an $S^1$-central
extension can be readily constructed by the methods of Poisson geometry. 
          
Applications of this construction to momentum map theory
and twisted K-theory will be discussed  elsewhere.

\section{Symplectic groupoids of  affine Poisson manifolds}
\subsection{General Construction}
Let $\frakg$ be a (finite or infinite dimensional) Lie algebra over
$\rr$,  and $\lambda \in \wedge^2\frakg^*$ a  Lie algebra
2-cocycle.
 Let $\widetilde{\frakg}=\frakg \oplus \reals$ be
the corresponding central extension.
 Assume that $S^{1}\lon \tilde{G}\stackrel{\pi}{\lon} G$
is  a central extension on the level of Lie groups, which 
exists if   $\omega_G \in \Omega^2 (G)$, the
left invariant closed two-form corresponding to $\lambda$,
is of integer class. 
It is well-known that  the transformation groupoid
$ \widetilde{\gm } (=\widetilde{G} \times \widetilde{\frakg}^* )\toto 
\widetilde{\frakg}^* $, where $\widetilde{G} $ acts on $\widetilde{\frakg}^*$
by coadjoint action:  $\tilde{g}\cdot \tilde{\xi}=
Ad_{\tilde{g}^{-1}}^* \tilde{\xi}, \ \forall \tilde{g}
\in \widetilde{G}, \xi \in \widetilde{\frakg}^*$
($Ad_{\tilde{g}^{-1}}^* $ stands for  the dual of $Ad_{\tilde{g}^{-1}}$),
 is a symplectic groupoid. 
The symplectic structure  on $\widetilde{G} \times \widetilde{\frakg}^* $
is the canonical cotangent symplectic structure when
$\widetilde{G} \times \widetilde{\frakg}^* $ is being    identified 
with $T^* \tilde{G}$ via the right translation.

Denote by $\chi : G\lon \frakg^*$ the group 1-cocycle
integrating the Lie algebra 1-cocycle $\lambda^{b}:\frakg\to \frakg^*$, 
$\la \lambda^{b}(v) , u\ra=\lambda (v, u), \ \forall v, u\in  \frakg$,
where  $G$ acts on $\frakg^*$ by  the coadjoint action.
We assume that $\chi$ exists, which is true, for instance,
when $G$ is simply connected.
Since $\tilde{G}$ is a central extension of $G$, its adjoint
action on $ \widetilde{\frakg}$ descends to  an
action of $ G$ given by
$g\cdot (X, t)= (Ad_{g}X,   \ t +\la \chi  (g^{-1}) , X\ra ),  \ \forall
g\in G, \ (X, t)\in \widetilde{\frakg}  \ (\cong \frakg\oplus \rr )$,
and therefore the induced coadjoint action is 
\begin{equation} 
\label{eq:coadjoint}
g\cdot (\xi , t)= (Ad_{g^{-1}}^*\xi + t \chi (g), \  t ), \ \ \forall 
g\in G, \ (\xi , t)\in \widetilde{\frakg}^* (\cong \frakg^* \oplus \rr ).
\end{equation} 

Embed $\frakg^*$ as a hyperplane  of $\widetilde{\frakg}^*$
via the map $\phi: \xi \lon (\xi, 1), \ \forall \xi \in \frakg^*$.
 Clearly $\frakg^*$ is a Poisson submanifold
of $\widetilde{\frakg}^*$ with the affine Poisson relation:
$\{l_{X}, l_{Y}\}=l_{[X, Y]}+ \lambda (X, Y),  \ \forall
X, Y\in \frakg $.
By Eq. \eqref{eq:coadjoint}, this hyperplane is invariant
under the coadjoint action of $G$, on which  it takes the form:
\begin{equation}
\label{eq:coad2}
g\cdot \xi = Ad_{g^{-1}}^*\xi + \chi (g), \ \ \forall g\in G, 
\xi \in \frakg^*.
\end{equation}  
Let $\gm $ be the corresponding transformation
groupoid $G\times \frakg^* \toto \frakg^* $.
 One may also form the transformation groupoid
$R: \tilde{G}\times \frakg^* \toto \frakg^* $ 
(here the $\tilde{G}$-action  on $\frakg^*$ is essentially 
the  same  action \eqref{eq:coad2}   composing with the group morphism
 $\pi : \tilde{G}\to G$).
Then $R\toto \frakg^* $ is   a subgroupoid of the symplectic
groupoid $\tilde{\gm} \toto \widetilde{\frakg}^*$ under the
natural embedding $i: (\tilde{g}, \xi )\to (\tilde{g}, \phi (\xi ))$. 
By $\theta_{\tilde{\gm}}$ and $\omega_{\tilde{\gm}} $,
 we denote the Liouville one-form and
symplectic two-form on  $\tilde{\gm} \  (\cong T^* \tilde{G} )$ respectively,
and set $\omega_R = i^*\omega_{\tilde{\gm}} \in
\Omega^2 (R ) $ and $\theta_R =i^* \theta_{\tilde{\gm}} \in \Omega^1 (R)$.
By $\tilde{\pi} : R\to \gm $, we denote the natural
projection: $\tilde{\pi} ( \tilde{g}, \xi ) =(\pi (\tilde{g} ), \xi )$,
$\forall ( \tilde{g}, \xi )\in \tilde{G}\times \frakg^* $.

\begin{them}
\label{thm:2.1}
\begin{enumerate}
\item The closed two-form  $\omega_R $ is  basic for the $S^1$-bundle
   $R\to \gm$, and therefore descends
to a  closed two-form  $\omega_{\gm}$ on $\gm$, i.e.,
$\omega_{R}=\tilde{\pi}^*\omega_{\gm}$.
\item $\omega_{\gm}$ is  symplectic and
 compatible with the groupoid structure so that
it defines a symplectic groupoid on $\gm$.
This is the symplectic groupoid  of the affine
Poisson manifold  $\frakg^*$.
\item $\tilde{\pi} : R\lon \gm $ is a $S^1$-central extension
of
Lie  groupoids.
\item $\theta_R$ is a prequantization connection form  on the
$S^1$-bundle
$\tilde{\pi } : R\to \gm $  compatible with the groupoid
structure, i.e., $\partial \theta_R =0$ and $d \theta_R =
\tilde{\pi}^* \omega_{\gm}$.
\end{enumerate}
\end{them}

In other words, $\theta_R\in C^2_{DR}(R\lcom ) $ is a pseudo-connection
on $R$  in the sense of
\cite{BX} with $\omega_{\gm}\in Z^3_{DR}(\gm\lcom )$
 being its pseudo-curvature.
Since $\gm\toto \frakg^*$ is a transformation groupoid,
its De Rham cohomology is the equivariant cohomology $H^*_{G}(\frakg^* )$,
 where $G$ acts
 on $\frakg^* $ by the gauge action \eqref{eq:coad2}.   
The Lie groupoid $S^1$-central extension $R\to \gm$ can be considered
as a geometrical model realizing  the class 
$[\omega_\gm]\in H^3_{G}(\frakg^* )$. Indeed, if $\RR$ and $\XX$
are the differential stacks corresponding to the Lie groupoids
$R\toto \frakg^*$ and $\gm  \toto \frakg^*$, respectively,
then $\RR$ is an $S^1$-gerbe over $\XX$ whose Dixmier-Douady
class is equal to $[\omega_\gm]\in H^3_{DR}(\XX )$ \cite{BX, BX2}.

\begin{numrmk}
{\em Note that $R\to \gm$ is indeed the pull back
$S^1$-central extension of $\tilde{G}\stackrel{\pi}{\lon} G$
via the groupoid morphism  $\psi$ from $\gm  (=G\times \frakg^* )
\toto \frakg^*$
to $G\toto \cdot$ defined by the natural  projection.
As a consequence, $\omega_\gm$ and $\psi^* \omega_G$ define 
the same class in $H^3_{G}(\frakg^* )$:  the 
Dixmier-Douady class of the  $S^1$-gerbe  $\RR\to \XX$.
It would be interesting to investigate whether
this class is  non-trivial when 
 $\lambda \in \wedge^2\frakg^*$ is assumed to be
a non-trivial  2-cocycle (otherwise it is obvious
that $\psi^* \omega$ is a trivial class). In the case of  loop groups
below, one indeed obtains a non-trivial class.}
\end{numrmk}

\subsection{Loop Group Case}
We will apply the above construction to the case
of loop groups.

Let $(\cdot, \cdot )$ be an ad-invariant  non-degenerate symmetric bilinear
form on $\frakg$. It is well-known that 
$(\cdot, \cdot )$ induces a Lie algebra 2-cocycle on the
loop Lie algebra $\lambda \in \wedge^2 (L\frakg^* )$ defined by \cite{PS}:
\begin{equation}
\label{eq:lambda}
\lambda (X , Y)=\frac{1}{2\pi}\int_{0}^{2\pi}(X(s), Y'(s))ds, \  \forall
X(s) , Y(s)\in L\frakg.
\end{equation}

By $\widetilde{L\frakg}$ we denote its corresponding  Lie algebra
central extension.
Assume that $\lambda $ satisfies  the integrability  condition (i.e.,
the corresponding closed two-form $\omega_{LG} \in  \Omega^2 (LG)^{LG}$
is of integer class). It defines  a loop group
central extension $S^1 \lon \widetilde{LG} \stackrel{\pi}{\lon} LG$.
By identifying $L\frakg^*$ with $L\frakg$ via the bilinear
form $(\cdot, \cdot )$, the 1-cocycle $\chi$ admits the form:
$\chi (g(s))=g' (s) g(s)^{-1} , \  \forall g(s)\in LG$,
and the gauge action  \eqref{eq:coad2} becomes
\begin{equation}
\label{eq:gauge}
g\cdot \xi = Ad_{g^{-1}}^*\xi +  g'  g^{-1}, \ \ \forall g\in LG, \ \xi\in
L\frakg.
\end{equation}
This is  the standard  gauge transformation when
$L\frakg$  is  identified with the space of connections on  the
trivial bundle over the 
unit circle $S^1$.

As above, we can form  the transformation
groupoids $\gm : LG \times L\frakg \toto  L\frakg $ and 
$R: \widetilde{LG} \times L\frakg \toto  L\frakg $, and
define $\omega_R\in \Omega^2 (R)$ and $\theta_R \in \Omega^1 (R)$.

According to  Theorem \ref{thm:2.1}, we  see
that the  closed two-form  $\omega_R $ is  basic   and  descends
to a  closed two-form  $\omega_{LG \times L\frakg}$ on $\gm$.

\begin{cor}
\label{thm:LG}
 $(LG \times L\frakg \toto L\frakg , \  \omega_{LG \times L\frakg} )$
 is a symplectic groupoid integrating the affine
Poisson structure on $L\frakg$. Moreover,  $\tilde{\pi}:
\widetilde{LG} \times L\frakg \to {LG} \times L\frakg$ is
 a $S^1$-central extension of Lie  groupoids,  on which
$\theta_R\in C^2_{DR}(R\lcom) $ defines a pseudo-connection
with $\omega_{LG \times L\frakg}$ being its pseudo-curvature.
\end{cor}


\section{An $S^1$-gerbe over $G/G$}
\subsection{AMM-groupoids}
\label{sec:AMM}
Let $G$ be a Lie group equipped with an ad-invariant non-degenerate symmetric
bilinear form $(\cdot , \cdot )$. Consider the transformation groupoid
${G\times G}\toto {G}$, where $G$ acts on itself by  conjugation.
As in \cite{AMM},  we denote by $\theta $ and $\bar{\theta}$
the left and right Maurer-Cartan forms on $G$ respectively, i.e.,
$\theta =g^{-1}dg$ and $\bar{\theta}=dg g^{-1}$.
Let $\Omega \in \Omega^{3}(G)$ denote the bi-invariant
3-form on $G$ corresponding to the Lie algebra 3-cocycle
$\frac{1}{12}(\cdot , [\cdot , \cdot ])\in \wedge^3 \frakg^*$,
and  $\omega \in \Omega^2 (G\times G )$  the two-form:

\begin{equation}
\label{eq:quasi}
\omega =-\half [(Ad_{x} g^* \theta , g^* \theta )
+(g^* \theta , x^{*}(\theta +\bar{\theta} ))],
\end{equation}
where $(g, x)$ denotes the coordinate in $G\times G$, and
$g^* \theta$ and $x^* \theta$ are, respectively,
  the $\frakg$-valued one-forms
on $G\times G$ obtained by pulling back
$\theta$ via the first and second projections, and similarly
for $x^* \bar{\theta} $.

A simple computation leads to 

\begin{prop}
\label{pro:AMM}
$\omega+\Omega $ is a  3-cocycle  of the De-Rham total
complex of the transformation groupoid ${G\times G}\toto {G}$,
and therefore it defines a class in the equivariant cohomology
$H^{3}_{G}(G)$.
\end{prop}

\begin{numrmk}
{\em \begin{enumerate}
\item When $G$ is a   compact simple Lie group with the  basic
form $(\cdot , \cdot )$,
$[\omega+\Omega ]$  
is a generator of   $H^{3}_{G}(G)$.
 In Cartan model, it corresponds
to the class defined by the $d_{G}$-closed equivariant 3-form
$\chi_{G}(\xi )=\Omega-\half \la \theta +\bar{\theta}, \xi \ra :
\ \frakg^* \lon \Omega^* (G), \ \forall \xi\in \frakg^*$.
\item In general, given a transformation groupoid 
${G\times M} \toto {M}$ ($G$ is assumed to be compact), and
a $d_{G}$-closed equivariant 3-form
$\chi_{G}=\Omega +E(\xi )$, where $\Omega \in
\Omega^{3}(M)$ is an invariant closed 3-form
on $M$ and $E: \frakg^* \lon \Omega^{1}(M)$
a $G$-equivariant linear map, 
 an explicit formula for a two-form $\omega \in \Omega^2 (G\times M )$ 
can be found \cite{Mein}, using the Bott-Shulman construction,
such that $\omega +\Omega \in Z^3_{DR}(G\times M \toto M)$
defines the same class of $\chi_{G}$. 
\end{enumerate}
}
\end{numrmk}

\subsection{$S^1$-central extensions}
Next we want to construct an $S^1$-central extension of Lie
groupoids
which  realizes the class of the 3-cocycle $\omega +\Omega$ 
 in   Proposition \ref{pro:AMM} as its Dixmier-Douady class.
 Since $\Omega\in \Omega^3 (G)$
is not exact, first of all we need to pass to
a Morita equivariant groupoid \cite{BX}. There are many different 
choices for such a groupoid.
 Basically, one needs to choose a surjective submersion
$f: M'\to G$ such that the pullback three form $f^* \Omega$ is exact
on $M'$. Then $\gm' : M'\times_{G, \alpha}\gm\times_{\beta, G}M' 
 \toto M'$ becomes  a groupoid, and the natural projection  from
$\gm' \toto M'$ to $\gm \toto G$
is  a Morita morphism \cite{BX}.
 For instance, one choice is to
take a good open cover $\{U_i \}$ of $G$.
Another choice, which is the one that we will pursue
in this Note, is the infinite dimensional manifold $L\frakg$ while 
$f$ is the holonomy map $\hol: L\frakg\lon G$, i.e., the time-$1$  map
of the differential equation:  $\hol_s (X)^{-1} \frac{\partial}{\partial s}
\hol_s (X) =X, \ \hol_0 (X)=e$.
  Then we have
$\hol^* \Omega  =d\mu$, where $\mu $ is the two-from
on $L\frakg$: $\mu=\half \int_0^1
( \hol_s^* \btheta , \frac{\partial}{\partial s} \hol_s^* \btheta )ds$ \cite{AMM}.

\begin{prop}
\begin{enumerate}
\item We have a Morita morphism $f$
 of Lie groupoids from
$LG \times L\frakg \toto L\frakg$ to
$G\times G\toto G$,
which is given by  $f (g, r)=(g(0), \hol (r) )$ on the space of
morphisms  and  by $f (r )=\hol (r)$ on the space of
objects,  $\forall g\in
LG, \ r\in L\frakg $.
\item Under the induced isomorphism $f\upst :H^3_{G} (G)
\longiso H^3_{LG}(L\frakg )$, $[\omega +\Omega ]$ goes to
$[\omega_{LG \times L\frakg}] $. Indeed we have
$$\omega_{LG \times L\frakg}-f\upst (\omega +\Omega )=\delta \mu.$$  
\end{enumerate}     
\end{prop}

As a consequence, we have

\begin{them}
Let $G$ be a Lie group equipped with an ad-invariant non-degenerate symmetric
bilinear form $(\cdot , \cdot )$. Assume that
$\lambda\in \wedge^2 (L\frakg^* )$ as in Eq. \eqref{eq:lambda}
satisfies the integrability condition.
Then the 3-cocycle $\omega +\Omega $  corresponds to an integer class
in $H^3_{G}(G)$.
The  $S^1$-central extension  $\widetilde{LG} \times L\frakg \toto L\frakg $
of ${LG} \times L\frakg \toto L\frakg $
corresponds  to an $S^1$-gerbe over the stack $G/G $ with the Dixmier-Douady
class $[\omega +\Omega ]\in H^3_{G}(G)$.                    
\end{them}

\begin{numrmk}
{\em If  there is a Morita morphism from  Lie 
 groupoid $\Gamma'\toto M'$ to $\Gamma\toto M$, then
these two groupoids are also  Morita equivalent in the sense
of \cite{Xu:1991sg}, which means  that
there is a bimodule.
Indeed, these two notions of Morita equivalence
are equivalent \cite{BX2}.
Morita equivalence via  bimodules is particularly
useful in constructing $S^1$-central extensions.
 It  allows one  to construct the $S^1$-central extension  of 
  one groupoid
in terms of an $S^1$-central extension of the other together
with a prequantization of the  bimodule. 
See \cite{BX2} for the details.
  This, for instance,
will lead to  a construction of an $S^1$-central extension
of the  Morita equivalent Lie groupoid $\gm'\toto M'$
when $M'=\cup U_{i}$ is an open covering as in \cite{Brylinski, Mein}
.}
\end{numrmk}

We end the paper with
the following proposition which  explicitly describes
the  equivalence bimodule between
the groupoids  ${G\times G}\toto G$ and
$LG\times L\frakg \toto { L\frakg}$.

\begin{prop}
The groupoids ${G\times G}\toto {G}$
and $LG\times L\frakg \toto { L\frakg}$ are
Morita equivalent 
in the sense of Definition 2.1 in \cite{Xu:1991sg},
where the bimodule $X$ can be  taken
as $G\times L\frakg$, and  $\rho : X\lon G$ and $\sigma :X
\lon L\frakg$ are given, respectively by
$\rho (g, r)=g\hol (r) g^{-1}$ and $\sigma (g, r)=r$,
$\forall (g, r)\in G\times L\frakg$.
 
The groupoid  ${G\times G}\toto {G}$ acts on
$X$ from the left by:
$(g_{1}, g_{2}) \cdot (g, r)= (g_{1}g , r), \  \forall (g_1 , g_{2})\in
G\times G, \ \ (g, r)\in G\times L\frakg$,
such that $g_2 =g\hol (r) g^{-1}$; while
$LG\times L\frakg \toto { L\frakg}$ 
acts from right:
$(g, r) \cdot (g(s), r')= (g g(0), r'), \ \forall (g, r)
\in G\times L\frakg, \ (g(s), r')\in LG\times L\frakg$,
such that $r=g(s)\cdot  r'$.
\end{prop}

\Acknowledgements{ We thank Eckhard Meinrenken and Alan
Weinstein for useful discussions, and Pierre Bieliavsky for his help
in preparing the French abbreviated version.
  We also  thank 
the Ecole Polytechnique 
 for hospitality and support of the research
summarized in this Note.   
This research is partially supported by   NSF
          grant DMS00-72171 and NSERC grant 22R81946.     }

\end{document}